\input amstex 
\documentstyle{amsppt}
\input bull-ppt
\keyedby{bull410/lic}

\topmatter
\cvol{29}
\cvolyear{1993}
\cmonth{July}
\cyear{1993}
\cvolno{1}
\cpgs{94-97}
\title A Counterexample to the\\ Rigidity Conjecture for 
Rings \endtitle
\author Raymond C. Heitmann\endauthor
\shortauthor{R. C. Heitmann}
\shorttitle{A counterexample to the Rigidity Conjecture 
for Rings}
\address Department of Mathematics, University of Texas,
Austin, Texas 78712\endaddress
\ml heitmann\@math.utexas.edu\endml
\date September 25, 1992\enddate
\abstract An example is constructed of a local ring and a 
module of
finite type and finite projective dimension over that ring 
such that
the module is not rigid. This shows that the rigidity 
conjecture is 
false.\endabstract
\endtopmatter

\document
\heading 1. Introduction\endheading
For a ring $R$, an $R$-module $M$ is said to be rigid if 
for every
$R$-module $N$ of finite type, whenever 
$\roman{Tor}_i^R(M,N)=0$, then
$\roman{Tor}_j^R(M,N)=0$ for all $j\ge i$. The concept was 
developed
in 1961 by Auslander, who showed that all modules over 
unramified regular
local rings were rigid \cite1 and who also demonstrated a 
number of other
results for modules over unramified regular local rings, 
all of which 
follow directly from the rigidity property and so would 
hold more
generally if a larger class of modules were known to be 
rigid. The
most notable of these came to be known as the zero divisor 
conjecture:
Let $R$ be a local ring, and let $M$ be an $R$-module of 
finite type and
finite projective dimension. If $x\in R$ is not a zero 
divisor on $M$,
then $x$ is not a zero divisor on $R$. (Auslander's 
original result
used regular sequences in place of nonzero divisors, but 
this more
elementary formulation, which appears in \cite{4, p.\ 8}, 
is equivalent.)

Auslander also noted the importance of finite projective 
dimension,
giving an example of a nonrigid module of infinite 
projective dimension;
thus, it was natural to focus on regular local rings, 
which have the
property that all modules have finite projective 
dimension. Auslander's
original question was answered in 1966 when Lichtenbaum 
proved the
rigidity conjecture for all regular local rings \cite5. 
However, because
of the numerous consequences of rigidity, it seemed 
natural to hope that
the property held for modules of finite projective 
dimension. Thus,
Peskine and Szpiro gave the conjecture its present form in 
\cite6:
If $R$ is a Noetherian ring and $M$ is an $R$-module of 
finite type and
finite projective dimension, then $M$ is rigid. In that 
article they offered
conjectures generalizing many of Auslander's results plus 
many other related
conjectures. The conjectures are primarily local in 
nature, as, for example,
a module is rigid if and only if it is locally rigid. For 
equicharacteristic
local rings of characteristic $p$ (both the ring and its 
residue field have
characteristic $p$\<), Peskine and Szpiro were successful 
in proving a
number of the conjectures, including the zero divisor 
conjecture. Hochster
extended this work to equicharacteristic zero, and his 
CBMS monograph
\cite4 gives a thorough presentation of the entire 
subject. There was,
however, no real progress on rigidity.

More recently, P. Roberts has been able to demonstrate 
some of the
conjectures which had previously eluded proof in the case 
of local
rings of mixed characteristic \cite8. In particular, he 
has shown that
Auslander's zero divisor conjecture and several other 
consequences of
rigidity are valid for all local rings. Thus, while the 
rigidity conjecture
has turned out to be false, many of its implications are 
true.

It is interesting to note that the example given here is 
not pathological.
It is an affine domain over a field constructed in a 
generic fashion.
It is not a complete intersection, and the conjecture thus 
is not resolved
for those rings. Hochster has suggested to me privately 
that he does not
believe this additional hypothesis will be helpful.

An alternate approach is to view rigidity as a property of 
the pair
$M,N$ rather than just $M$. For regular local rings, 
necessarily both $M$
and $N$ have finite projective dimension, and so, in this 
context, it
seems appropriate to assume that both modules have finite 
projective
dimension. This question likewise remains open.

\heading 2. The example\endheading
To settle this conjecture in the negative, we will 
construct a Noetherian ring
$R$ together with $R$-modules $M$ and $N$ such that 
$pd_RM=2$, $N$ is a
module of finite length, $\roman{Tor}_1^R (M,N)=0$, and 
$\roman{Tor}_2^R
(M,N)\ne 0$. Lemma 1 is not really new, and to avoid 
disturbing the flow
of the construction, we will defer its proof to the end of 
the article.

\thm{Lemma 1} There exists an affine $K$-algebra $R$ where 
$K$ is a field,
matrices $X=[x_{ij}]$ and $Y=[y_{ij}]$ of dimensions 
$2\times 4$ and
$4\times 8$ respectively, and a maximal ideal $P$ of $R$ 
such that\,\RM:
\roster
\item $0\to R^2 @>X>> R^4 @>Y>> R^8$ is exact\,\RM;
\item $\{\overline x_{ij},\overline y_{ij}\}$ is a 
linearly independent
subset of $P/P^2$\RM; and
\item $R/P=K$. 
\endroster
\ethm

\dfn{Definition} Let $S=K[s,t]/(s^2,st,t^2)$. Also let
$N=S^2/((t,0),(0,s), 
(s,t))S$, a length 3 $S$-module. Note that $(s,t)N$ is a 
length 1 submodule.
\enddfn

\thm{Lemma 2} There is a homomorphism $R\to S$ which takes
$$
X\ \text{to}\ \overline X=\pmatrix s& 0 & t & 0\\
0 & s & 0 & t\endpmatrix\quad\text{and}\quad Y\ \text{to}\ 
\overline Y
=\pmatrix s& 0 & 0 & 0 & t & 0 & 0 & 0\\
0 & s & 0 & 0 & 0 & t & 0 & 0\\
0 & 0 & s & 0 & 0 & 0 & t & 0\\
0 & 0 & 0 & s & 0 & 0 & 0 & t\endpmatrix.
$$
This map gives $N$ an $R$-module structure.\ethm

\demo{Proof} It suffices to find a homomorphism $R/P^2\to 
S$. 
So we may assume $R=K\oplus P$ where $P$ is a $K$-vector 
space
and $S=K\oplus(s,t) K$. Then a $K$-algebra homomorphism is 
determined by a linear transformation $P\to (s,t)K$. Now, as
$\{\overline x_{ij},\overline y_{ij}\}$ is linearly 
independent,
there is clearly such a transformation with the desired
$\overline X,\overline Y$.\qed
\enddemo

\thm{Theorem 3} With the above notation and $M=R^8/(R^4)Y$,
$pd_RM=2$ and $\roman{Tor}_1^R (M,N)=0$, but 
$\roman{Tor}_2^R
(M,N)\ne 0$.\ethm

\demo{Proof} Since $0\to R^2 @>X>> R^4 @>Y>> R^8 \to M$ is a
projective resolution of $M$, it follows that $pd_RM\le 
2$. As
each entry of $X$ is in $P$, the injection does not split, 
and
so we actually have equality.

To compute Tor, note that $_-\otimes_R 
N=(_-\otimes_RS)\otimes_SN$.
Thus we can compute $\roman{Tor}_i^R(M,N)$ as the homology 
modules
of the complex of $S$-modules $0\to N^2 @>\overline 
X\otimes N>>
N^4 @>\overline Y\otimes N>> N^8$. Then it is easy to see 
that the
image of $\overline X\otimes N$ is $(s,t)N^4$ and the 
image of 
$\overline Y\otimes N$ is $(s,t)N^8$. As the length of 
$N^4$ is
12, this forces the sequence to be exact at $N^4$, and so 
$\roman{Tor}
_1^R(M,N)=0$. On the other hand, 
$\roman{Tor}_2^R(M,N)=(s,t)N^2\ne 0$.\qed
\enddemo

\demo{Proof of Lemma \rm1} Begin with generic matrices 
$X,Y$ of dimensions
$2\times 4$ and $4\times 8$ respectively. Let $K$ be a 
field, let
$R_1=K[x_{ij},y_{ij}]$, and let $J$ be the ideal of $R_1$ 
generated by
the entries of $XY$ and all $3\times 3$ minors of $Y$. Let 
$R_2=R_1/J$.
This is the ring which Bruns calls $K((2,4,8),(2,2))$ 
\cite{2, pp.\ 53--55}.
Let $f$ be the $2\times 2$ minor of $X$ determined by 
columns 3 and 4, 
and let $\{g_\sigma\}$ be the set of all $2\times 2$ 
minors of $Y$ which
use the first two rows (so there is one $g_\sigma$ for 
each pair of columns).
Let $R=R_2[\{g_\sigma/f\,\}]$. It is clear from \cite2 
that $0\to R^2
@>X>> R^4 @>Y>> R^8$ is exact. (The elements $g_\sigma/f$ 
are the
$u$-elements in Bruns's paper.)

Next, we claim $R$ is a graded ring. Let $\deg (x_{ij})=2$ 
for each
$i,j$, and let $\deg(y_{ij})=3$ and 
$\deg(g_\sigma/f\,)=2$. Then it
is clear that $R$ is a graded polynomial ring modulo a 
homogeneous 
ideal and so is graded. Let $P$ be the elements of 
positive degree.
Then $R/P=K$ and $P$ is maximal. Finally, all of the 
relations between
the generators have degree at least 4, and so $\{\overline 
x_{ij},
\overline y_{ij}\}$ is a linearly independent subset of 
$P/P^2$. \qed
\enddemo

\rem{Remarks} The salient features of the example are 
unchanged if we
localize at $P$ or localize and complete.

The Betti numbers (ranks of the free modules in the 
projective 
resolution) of $M$ are $\langle 8\ 4\ 2\rangle$. A similar 
example 
can be constructed with Betti numbers $\langle 9\ 3\ 
1\rangle$;
in that case, $N$ will be a module of length 4.

In \cite{3, p.\ 287} it was conjectured that rigidity was 
false and a
counterexample might exist with Betti numbers $\langle b\ 
b\ 1\rangle$
for $b>2$. We do not rule this out. However, such an 
example will
necessarily be more complex. The length of $N$ must be 
considerably
larger, and $N$ will not be annihilated by $P^2$ or even 
$P^3$.
\endrem

\heading Acknowledgment\endheading
The author thanks David Eisenbud, Mel Hochster, and Craig 
Huneke for
their assistance in the preparation of this article.

\enddocument